\documentclass[10pt]{amsart}
\usepackage{amsmath,amscd,amssymb,epsfig}
\usepackage{amsmath,amscd,amssymb,verbatim,paralist}

\usepackage{amsmath}
\usepackage{amsthm}
\usepackage{amssymb}
\usepackage{verbatim}
\usepackage{subfigure}
\usepackage{verbatim}
\usepackage{url}
\numberwithin{equation}{section}

\def\displayandname#1{\rlap{$\displaystyle\csname #1\endcsname$}%
                      \qquad \texttt{\char92 #1}}

\makeatletter
\def\url@leostyle{%
  \@ifundefined{selectfont}{\def\UrlFont{\sf}}{\def\UrlFont{\small\ttfamily}}}
\makeatother
\newcommand{\silsul}[2]{\genfrac\{\}{0pt}{}{#1}{#2}}
\newcommand{\ribua}[2]{\genfrac[]{0pt}{}{#1}{#2}}

\newtheorem{thm}{Theorem}[section]
\newtheorem{pro}[thm]{Proposition}
\newtheorem{lem}[thm]{Lemma}
\newtheorem{cla}[thm]{Claim}

\newtheorem{cor}[thm]{Corollary}
\newtheorem{fac}[thm]{Fact}
\newtheorem{obs}[thm]{Observation}

\newtheorem{df}[thm]{Definition}

\newtheorem{nt}[thm]{Note}
\newtheorem{nota}[thm]{Notation}

\begin{document}

\bibliographystyle{acm}

%\title{Generator Sets for the Alternating Group}
%
%\author{Aviv Rotbart \thanks{This paper is part of the author's Master Thesis, which was written at Bar-Ilan University under the supervision of R. M. Adin and Y. Roichman.}\\
%Department of Mathematics\\
%Bar-Ilan University\\
%Ramat Gan 52900, Israel\\
%{\em aviv.rotbart@gmail.com}} 
%\maketitle

%\tableofcontents

\begin{center}
\LARGE{Generator Sets for the Alternating Group}

\large{Aviv Rotbart \footnote{Department of Mathematics, Bar-Ilan University, Ramat Gan 52900, Israel.
\\ Email:  \tt{aviv.rotbart@gmail.com}} 
\footnote{This paper is part of the author's MSc thesis, which was written at Bar-Ilan University under the supervision of Prof. R.\ M.\ Adin and Prof. Y.\ Roichman.} 
} 
\end{center}

\begin{center}
\section*{Abstract}
\end{center}

%\begin{abstract}

{\small Although the alternating group is an index 2 subgroup of the
symmetric group, there is no generating set that gives a Coxeter structure
on it.  Various generating sets were suggested and studied by Bourbaki,
Mitsuhashi, Regev-Roichman, Vershik-Vserminov and others.  In a recent
work of Brenti-Reiner-Roichman it is explained that palindromes in
Mitsuhashi's generating set play a role similar to that of reflections in
a Coxeter system.

We study in detail the length function with respect to the set of
palindromes. Results include an explicit combinatorial description, a
generating function, and an interesting connection to Broder's restricted
Stirling numbers.}

%\end{abstract}

%\maketitle

\section{Introduction}
The study of parameters (statistics) of the symmetric group and other related groups is a very active branch of combinatorics in recent years. A major step was made about one hundred years ago, when MacMahon~\cite{MM} showed that the parameters \emph{major index} and \emph{inversion number} are equi-distributed on the symmetric group, $S_n$. This important result is the foundation of the field, and stimulated many subsequent generalizations and refinements. \\

It is well known that statistics on $S_n$ may be defined via its Coxeter generators (simple reflections) $\{ s_i=(i,i+1)\mid 1\le i \le n-1\}$, or via the transpositions (reflections) $\{ t_{ij}=(i,j)\mid 1\le i<j \le n\}$. Unfortunately, the alternating group $A_n \subseteq S_n$ is not a Coxeter group. Our goal is to study generating sets for the alternating group that play a role similar to that of reflections in the symmetric group, and to explore the combinatorial properties of $A_n$ based on these sets.\\

A good candidate is the set $\{s_1s_{i+1}=(1,2)(i+1,i+2)\mid 1< i < n-1\}$. Mitsuhashi~\cite{MI} pointed out that these generators for the alternating group play a role similar to that of the above Coxeter generators of $S_n$. Regev and Roichman~\cite{RR} describe a canonical presentation of the elements in $A_n$ based on this set. They also calculate the generating functions of length and other statistics, with respect to this set of generators.\\

Our work deals with $A_n$-statistics calculated with respect to a new set of generators, $\{s_1t_{ij}=(1,2)(i,j)\mid 1\le i<j \le n\}$. This set consists of palindromes in Mitsuhashi's generators discussed above. As Brenti, Reiner and Roichman [3] explain, these palindromes play a role similar to that of reflections in the symmetric group. The following diagram describes the relations between the four generating sets mentioned above. 

$$
\begin{matrix}
  \text{$S_n$-Coxeter} & \underrightarrow{\text{conjugate by $S_n$}} & \text{$S_n$-Transpositions} \\
  C=\{ s_i\mid 1\le i \le n-1\}&          & T=\{ t_{ij}\mid 1\le i<j \le n\}&\\ \\
 \Biggm\downarrow{\text {multiply by $s_1$}}&                   &\Biggm\downarrow{\text {multiply by $s_1$}} & \\ \\
  \text{$A_n$-Coxeter (Mitsuhashi)}&  \underrightarrow{\text{take palindromes}} & \text{$A_n$-Transpositions} \\
    C(A_n)=\{ s_1s_{i+1}\mid 1< i \le n-1\}&         & T(A_n)=\{s_1t_{ij}\mid 1\le i<j \le n\}&\\
\end{matrix}
$$ \\

Various aspects of the generating set $T(A_n)$ are studied in this work, including: canonical forms of elements in $A_n$ with respect to $T(A_n)$; length of elements and the relation between length and number of cycles; a generating function for length expectation and variance for length; and finally a connection with Broder's restricted Stirling numbers~\cite{BRO}.

The methods used in this work include: manipulations on generating functions of Stirling numbers; theoremes on the number of cycles of a permutation and on the number of permutations of a given length in $S_n$ ; bijections between certain subsets of permutations in $A_n$, $S_n$ and permutations with Broder's property (see Definition \ref{dfrstr1}). \\

The paper is organized as follows : Detailed background and notations for the symmetric and alternating groups, as well as for Stirling numbers, is given in Section 2. In Section 3 we present the main results achieved in this work. The \emph{A} canonical presentation is analyzed in Section 4. In section 5 we discuss refined counts of permutations in $A_n$, while the relation between length and the number-of-cycles statistic is analyzed in Section 6. In Section 7 we calculate the generating function of length with respect to the generating set $T(A_n)$. The expectation and variance of the length function are studied in Section 8. The relation between our results and restricted Stirling numbers is analyzed in Section 9.

\section{Background}

\subsection{The Symmetric Group}
In this subsection we present the main notations, definitions and theorems on the symmetric group, denoted $S_n$.

\begin{nota}
Let $n$ be a nonnegative integer, then $[n]:=\{1,2,3,\dots ,n\}$ (where $[0]:=\emptyset$). 
\end{nota}

%\begin{df}
%Let $\mathbb{N}$ be the set of all natural numbers. A \emph{permutation of order $n\in \mathbb{N}$} is a bijection $v:[n]\to [n]$.
%\end{df}
%
%\begin{df}
%a \emph{cycle} is a permutation of order $n$ which maps the elements of some subset $S\subseteq [n]$ to each other in a cyclic fashion, while fixing (i.e., mapping to themselves) all other elements. The set $S$ is called the orbit of the cycle.
%\end{df}
%
%The following is a basic and very useful theorem in combinatorics.
%\begin{thm}
%Each permutation can be written as a product of disjoint cycles, unique up to the order of the cycles and a cyclic order in each cycle.
%\end{thm}
%
%\begin{df}
%Call the above presentation of a permutation $v$ the \emph{Cycle Notation} of $v$.
%\end{df}
%
%\begin{rem}
%Permutations are traditionally written in a two-line notation
%$$
%v=
%            \left(%
%                \begin{array}{ccccc}
%                    1 & 2 & 3 & \ldots & n \\
%                    v(1)& v(2) & v(3) & \ldots & v(n) \\
%                \end{array}
%            \right),
%$$
%however we will use the cycle notation which is more convenient.
%$$
%v=(v(1),v(2),v(3),\ldots,v(n)).
%$$
%TODO: maybe example... show also cycle notation.
%\end{rem}
%
%\begin{exa}
% $v=
%            \left(%
%                \begin{array}{cccccc}
%                    1 & 2 & 3 & 4 & 5 & 6 \\
%                    6 & 2 & 4 & 5 & 3 & 1 \\
%                \end{array}%
%            \right)$ will be written as $v=(16)(2)(345)=(16)(345)$. Cycles of length $1$ can be omitted in this notation.
%\end{exa}

\begin{df}
Denote by $\mathbb{N}$ the set of natural numbers.
The \emph{symmetric group} on $n\in \mathbb{N}$ letters (denoted $S_n$)  is the group consisting of all permutations on $n$ letters, with composition as the group operation.
\end{df}

\begin{df}
Given a permutation $v\in S_n$, we say that a pair $(i,j)\in [n]\times [n]$ is an \emph{inversion} of $v$ if $i<j$ and $v(i)>v(j)$. If $(i,i+1)$ is a transposition of $v$ then it is called an \emph{adjacent transposition}.
\end{df}

\begin{df}
The \emph{Coxeter generators} of $S_n$ are 
$$
\{ s_i=(i,i+1)\mid 1\le i \le n-1\},
$$ 
i.e., all the adjacent transpositions.
\end{df}

It is a well-known fact that the \emph{symmetric group} is a \emph{Coxeter
group} with respect to the above generating set. The following natural statistic describes the length of permutations in the symmetric group, with respect to the Coxeter generating set:

\begin{df}\label{lengthSnC} The \emph{length} of a permutation $v\in S_n$ with respect to the Coxeter generators is defined to be 
$$
\ell_{C}(v):=\min\{\ r \geq 0 \ | \  v=s_{i_1} \ldots s_{i_r} \mbox{ for some } i_1,\ldots,i_r \in [n-1] \ \}.
$$
\end{df}

\begin{df}
The \emph{inversion number} of $v\in S_n$ is 
$$
inv(v):=|\{(i,j) \ | \ 1 \leq i<j\leq n, \ v(i)>v(j)\}|
$$
\end{df} 

\begin{fac}
For each $v\in S_n$,
$$
inv(v)=\ell_{C}(v)
$$
\end{fac}

Another important set of generators for $S_n$ is the set of all transpositions.

\begin{nota}
Denote by $T$ the set of all transpositions in $S_n$, i.e. 
$$
T=\{(i,j)\mid 1\le i<j \le n\}
$$
\end{nota}

The definition of length with respect to $T$ is similar.

\begin{df}\label{lengthSnT}
Let $v\in S_n$, then  
$$
\ell_{T}(v):=\min\{\ r \geq 0 \ | \  v=t_1
\ldots t_r,\quad t_i \in T \ \}.
$$
\end{df}

A well known result describes the connection between the number of cycles and this length statistics in $S_n$.
\begin{thm}
If $cyc(v)$ is the number of cycles in $v\in S_n$, then
$$
\ell_{T}(v) + cyc(v)=n
$$
\end{thm} 

This result will be useful in some of the proofs in this work.
%\begin{df}
%The \emph{sign} of a permutation $v\in S_n$ is defined as follows.
%$$
%sign(v)=(-1)^{\ell(v)}
%$$
%$v$ is called an \emph{even permutation} in $S_n$ if its sign is $1$, and an \emph{odd permutation} otherwise.
%\end{df}
%
%\begin{nt}
%The sign statistic can be defined similarly with $\ell_{T}(\cdot)$ instead of $\ell(\cdot)$ and the two definitions are equivalent for each $v\in S_n$. 
%\end{nt}

\subsection{The Alternating Group}\label{altgrp_sec}
In this section we will define the alternating group, which is a subgroup of the symmetric group. We will also describe a known generating set for this group and the corresponding generating function of length.

\begin{df}
The \emph{Alternating Group} on $n$ letters, denoted $A_n$, is the group consisting of all even permutations in the symmetric group $S_n$; i.e., $A_n:=\{v\in S_n\mid sign(v)=1\}$. 
\end{df}

%\begin{fac}
%The alternating group $A_n$ is a normal subgroup of $S_n$; i.e., for all $w\in S_n$, $w A_n w^{-1}=A_n$
%\end{fac}
%
%\begin{fac}
%The cardinality of $A_n$ is 
%$$
%|A_n|=\frac{|S_n|}{2}=\frac{n!}{2}
%$$
%\end{fac}

% TODO:  cite
Following Mitsuhashi~\cite{MI} we let
$$
a_i:=s_1s_{i}=(1,2)(i,i+1)\qquad (2\le i\le n-1).
$$

The set $C(A_n):=\{a_i\ | \ 2\le i\le n-1\}$ generates the alternating group on $n$ letters, $A_{n}$.

\medskip

Regev and Roichman~\cite{RR} used Mitsuhashi's generators to describe a covering map  $f:A_{n+1} \to S_n$, which allows us to translate $S_n$-identities into corresponding $A_{n+1}$-identities. They gave a formula for the generating function of length with respect to these generators.

\begin{pro}~\cite[Thm.~6.1]{RR}
$$
\sum_{w\in A_{n+1}}q^{\ell_{C(A_n)}(w)} = (1+2q)(1+q+2q^2)\cdots (1+q+\ldots + q^{n-2}+2q^{n-1})
$$
Where $\ell_{C(A_n)}(\cdot)$ is the length with respect to Mitsuhashi's generators.
\end{pro}

\subsection{Stirling Numbers}
For basic properties of Striling numbers the reader is reffered to ~\cite{St}. In this subsection we will describe one important generalization of them, Broder's ~\cite{BRO} restricted Stirling numbers.

\begin{df}\label{dfrstr1}
The unsigned r-restricted Stirling number of the first kind, denoted $\ribua{n}{k}_r$, is the number of permutations of the set $\{1,2,...,n\}$ with $k$ disjoint cycles, with the restriction that the numbers $1, 2, ..., r$ belong to distinct cycles. The case $r=1$ gives the usual unsigned Stirling numbers of the first kind.
\end{df}

\begin{df}
The Kronecker delta function is defined as follows.
$$
\delta_{i,j}=
\begin{cases}
1,& \text{if } i=j; \\
0,& otherwise.
\end{cases}
$$
\end{df}

%clabla
\begin{cla}
r-Stirling numbers of the first kind satisfy the same recurrence relation as unsigned Stirling numbers of the first kind, except for the initial conditions:
$$
\ribua{n}{k}_r=(n-1){\ribua{n-1}{k}}_r+{\ribua{n-1}{k-1}}_r \qquad (r<k<n)
$$
with the following boundary conditions:
\begin{align*}
&\ribua{n}{k}_r=0, & (k<r \text{ or } n<k); \\
&\ribua{n}{r}_r=\frac{(n-1)!}{(r-1)!},  &(r \leq n); \\
&\ribua{n}{n}_r=1, &(r \leq n). \\
%&\ribua{n}{0}_r=0.  &(\forall n)
\end{align*}

%(n-1){\ribua{n-1}{k}}_r+{\ribua{n-1}{k-1}}_r& n>r

\end{cla}

%TODO:cite
\begin{thm}~\cite[\S 6.9]{BRO}
The generating function of unsigned r-restricted Stirling numbers of the first kind is
$$
\sum_{k=0}^{n} {\ribua{n}{k}}_r\cdot x^{k}=
\begin{cases}
x^r(x+r)(x+r+1)\cdots (x+n-1)& \text{if } 1\leq r\leq n; \\
0& otherwise.
\end{cases}
$$
\end{thm}
%Compare this result with Theorem \ref{eq8}.

\begin{df}\label{dfrstr2}
The r-restricted Stirling number of the second kind, denoted ${\silsul{n}{k}}_r$, is the number of ways to partition the set $\{1,2,...,n\}$ into $k$ nonempty disjoint subsets with the restriction that the numbers $1, 2, ..., r$ belong to distinct subsets. The case $r=1$ gives the usual Stirling numbers of the second kind.
\end{df}

\begin{cla}
r-restricted Stirling numbers of the second kind satisfy the same recurrence relation as Stirling numbers of the second kind, except for the initial conditions.
$$
\silsul{n}{k}_r={k\cdot \silsul{n-1}{k}}_r+{\silsul{n-1}{k-1}}_r \qquad (r<k<n)
$$
with the following boundary conditions:
\begin{align*}
&\silsul{n}{k}_r=0, & (k<r \text{ or } n<k); \\
&\silsul{n}{k}_r=r^{n-r},  &(r\leq n); \\
&\silsul{n}{n}_r=1, &(r\leq n). \\
%&\silsul{n}{0}_r=0.
\end{align*}

%$$
%{\silsul{n}{k}}_r=
%\begin{cases} 
%0& n<r, \\
%\delta_{k,r}&  n=r, \\
%k\cdot{\silsul{n-1}{k}}_r+{\silsul{n-1}{k-1}}_r& n>r
%\end{cases} 
%$$
\end{cla}

\begin{thm}~\cite[\S 6.10]{BRO}
The generating function of r-restricted Stirling numbers of the second kind is
$$
\sum_{n=0}^{\infty} {\silsul{n}{k}}_r\cdot x^{n}=
\begin{cases}
\frac{x^k}{(1-rx)(1-(r+1)x)\cdots (1-kx)},& \text{if } 1\leq r\leq k; \\
0,& otherwise.
\end{cases}
$$
\end{thm}
%Compare this result with Theorem \ref{eq7}.

Restricted Stirling numbers of the first and second kind satisfy the same orthogonality relation as the usual (unsigned) Stirling numbers, as described in the following theorem.

\begin{thm}\label{thmres}~\cite[\S 4.5]{BRO}
$$
\sum_{k=0}^{n} {\ribua{n}{k}}_r\cdot{\silsul{k}{m}}_r\cdot (-1)^k=
\begin{cases}
(-1)^n\cdot\delta_{m,n}, & \text{if } r\leq m\leq n; \\
0,& otherwise.
\end{cases}
$$
\end{thm}
\subsection{Harmonic Numbers}
\begin{df}\label{def_harmc}
The \emph{$n$-th harmonic number}, denoted by $H_n$, is the sum of the reciprocals of the first $n$ positive integers:
$$
H_n= 1+\frac{1}{2}+\frac{1}{3}+\cdots + \frac{1}{n}.
$$
\end{df}
\begin{df}\label{def_genharmc}
The \emph{generalized $n$-th harmonic number of order m}, denoted $H_{n,m}$, is 
$$
H_{n,m}= 1+\frac{1}{2^m}+\frac{1}{3^m}+\cdots + \frac{1}{n^m}.
$$
\end{df}

\section{Main Results}
In this section we present the main results of this paper. Details and proofs will be given in Sections 4-8.

%We begin section 4 with canonical presentation of elements in $A_n$ and prove some properties of this presentation.

Let
$$
a_{ij} := s_1t_{ij}=(12)(ij) \qquad (1 \le i < j \le n).
$$
The set of \emph{A-transpositions} 
$$
T(A_n):=\{a_{ij} \mid 1 \le i < j \le n \}
$$ 
generates the alternating group on $n$ letters. The length of an element $v\in A_n$ can be naturally defined with respect to the above generators:
$$
 \ell_{T(A_n)}(v)= \min \{k\geq 0 \mid v= v_1 \cdots v_k,\quad v_i\in T(A_n)\}
$$

%\begin{df} (Definition \ref{def_newgen})
%for $i \geq 3$ we define the following set of permutations
%\begin{equation*}
% R_i= \{ (12)(ji) \mid 1 \leq j<i \} \cup \{ e \} 
%\end{equation*}
%\end{df}
%
%\begin{thm} (Theorem \ref{thm_canon})
%For any $v \in A_n$, $n \geq 3 $, there exist elements $v_i \in R_i$, $3 \leq i \leq n$, such that $v=v_3\cdots v_n$, and this presentation is unique. Call it the canonical presentation of $v$. 
%\end{thm}
%Now we can define the length of a permutation with respect to the above canonical presentation.
%
%\begin{df} (Definition \ref{deflenm})
%If $v\in A_n$ has the canonical presentation $v=v_3 \cdots v_n$, let 
%$$
%\ell_m(v):=\#\{i\mid v_i\neq e\}
%$$
%\end{df}
%
%We then show that the definition of length according to the canonical presentation coincides with the natural definition of length with respect to the generators set $T(A_n)$ defined above (see Theorem \ref{thm3} below).
%
%$$
%\ell_m(v) = \ell_{T(A_n)}(v).
%$$

\begin{nota}
Denote by $a(n,m)$ the number of elements of length $m$ in $A_n$.
\end{nota}

Our first result is a Stirling-type recursion for $a(n,m)$.
\begin{pro} (Corollary \ref{conc_len})
$$
a(n,m)=(n-1)\cdot a(n-1,m-1) + a(n-1,m) \qquad (0<m<n) 
$$

with boundary conditions $a(n,0)=1$ for $n\geq 0$, and $a(n,n)=0$ for $n>0$. 
\end{pro}

%Section 6 discusses the relation between the \emph{length} and \emph{number of cycles} statistics. We prove that all the permutations in $A_n$ of length $k$ have the same number of cycles, and denote this number by $m(n,k)$. We then give an explicit formula for this number. A permutation $v\in A_n$ is of length $\ell_{T(A_n)}(v)=k$ if and only if the number of cycles in $v$ is
%$$
%n-k+\frac{(-1)^k-1}{2}=\\
%\begin{cases}
%n-k& \text{if $k$ is even}, \\
%n-k-1&  \text{if $k$ is odd}
%\end{cases} 
%$$
%
%\begin{thm} (Theorem \ref{thm_cyc})
%$$
%m(n,k)=n-k+\frac{(-1)^k-1}{2}=\\
%\begin{cases}
%n-k& \text{k is even}, \\
%n-k-1&  \text{k is odd}
%\end{cases} 
%$$
%\end{thm}
%
The following result relates the length function to the cycle number.

%\begin{df}
%For $v\in A_n$ given as product of disjoint cycles, let $cyc(v)$ be the number of cycles in $v$.
%\end{df}

\begin{pro} (Corollary \ref{conc_cyc})
Let $v\in A_n$, $n\geq 2$, Then
\begin{equation*}
\ell_{T(A_n)}(v)= \\
\begin{cases}
n-cyc(v)&  \text{if 1,2 are in different cycles of v;} \\
n-cyc(v)-1&  \text{if 1,2 in the same cycle of v.}
\end{cases}
\end{equation*}
\end{pro}
(For $n\leq 2$, $A_n$ contains only the identity permutation.)

Note that the length function $\ell_{T(A_n)}(v)$ is odd if and only if $1,2$ are in the same cycle in $v$ (see Corollary \ref{pro_nice}).

\begin{pro} (Theorem \ref{thm_genfunc}) For $n\geq 2$,
\begin{align*}
\sum_{v\in A_n} x^{\ell_{T(A_n)}(v)}
&=\sum_{k=0}^n a(n,k)\cdot x^{k} \\
&=(1+2x)(1+3x)\cdots (1+(n-1)x) \\
&= \prod_{t=2}^{n-1}(1+tx) \nonumber
\end{align*}
\end{pro}

\begin{thm} (Theorem \ref{thm_exp_var})
%Let
%\begin{align}\label{eq_genfunc}
%\sum_{v\in A_n} x^{\ell_{T(A_n)}(v)}
%= \prod_{t=2}^{n-1}(1+tx)   \nonumber
%\end{align}
%be the generating function of length in $A_n$ with respect to the generating set $T(A_n)$, then
The expected value of $\ell_{T(A_n)}$ is
$$
E[\ell_{T(A_n)}]=n-H_n-\frac{1}{2}
$$
and its variance is
$$
Var[\ell_{T(A_n)}]=H_n-H_{n,2}-\frac{1}{4}
$$
\end{thm}

Finally, we discuss a certain generalization of Stirling numbers and relate it to our statistic $a(n,k)$. The generalization discussed is Broder's restricted Stirling numbers~\cite{BRO}, see Definitions \ref{dfrstr1} and \ref{dfrstr2}. The connection was initially established using the On-Line Encyclopedia of Integer Sequences~\cite{EIS}.

\begin{pro} (Theorem \ref{thm_str}) For $0\leq k\leq n-2$,
$$
a(n,k)={\ribua{n}{n-k}}_2
$$
\end{pro}

%Hence the inverse matrix of signed $a(n,k)$ (or signed $\ribua{n}{k}_2$) can be calculated.
%
%\begin{cor} (Corollary \ref{con_rstr2}, see ~\cite[\S 4]{BRO})
%Let 
%\begin{align*}
%&A=[a_{i,j}]_{n\times n}, \qquad a_{i,j}=a(i,j)\cdot (-1)^{n-j} \\
%&B=[b_{i,j}]_{n\times n},\qquad   b_{i,j}=\silsul{i}{j}_2
%\end{align*}
%be the matrices of $a(n,k)$ values with alternating signs, and of $2$-restricted Stirling numbers of the second kind (see Definition \ref{dfrstr2}) respectively. Then 
%$$
%A\cdot B= I.
%$$

%inverse matrix of signed $a(n,k)$ (or signed ${\ribua{n}{k}}_2$) is exactly ${\silsul{n}{k}}_2$, i.e. 
%$$
%\Vert a(n,k)\cdot (-1)^k \Vert \times \Arrowvert \silsul {n}{n-k}_2 \Arrowvert = \Vert \delta_{i,j}\cdot (-1)^i \Vert 
%$$
%\end{cor}

\medskip

\section{The $A$ Canonical Presentation}
%\begin{rem} 
%In this section we present the main results obtained in this paper according to their order of appearance. 
%\end{rem}
In this section we consider a canonical presentation of elements in $A_n$ by the corresponding $s_1t_{ij}$ generators.

\subsection{A Generating Set for $A_n$}\label{subsdefan}
We let
$$
a_{ij} := s_1t_{ij}=(12)(ij) \qquad (1 \le i < j \le n).
$$

%The set of A-transpositions $T(A_n):=\{a_{ij} \mid 1 \le i < j \le n \}$ generates the alternating group on $n$ letters.
Denote by $T(A_n):=\{a_{ij} \mid 1 \le i < j \le n \}$ the set of A-transpositions.

\begin{df}\label{def_newgen}
For $n \geq 3$ define the following subset of permutations in $A_n$:
\begin{equation*}
 R_n= \{ (12)(in) \mid 1 \leq i<n \} \cup \{ e \}. 
\end{equation*}
\end{df}

%\begin{nt}
%$R_n$ is a subset of the generating set of the alternating group on $n$ letters.  $R_n=(T(A_n)\setminus T(A_{n-1}))\cup \{e\} .$
%\end{nt} 
\begin{nt}
$R_n$ is a subset of $T(A_n)$.  $R_n=(T(A_n)\setminus T(A_{n-1}))\cup \{e\} .$
\end{nt}

\begin{thm}\label{thm_canon}
Let $v \in A_n$, $n \geq 3 $. Then there exist unique elements $v_i \in R_i$, $3 \leq i \leq n$, such that $v=v_3\cdots v_n$. Call it the canonical presentation of $v$.
\end{thm}

%\begin{df}
%Call the above $v=v_3\cdots v_n$ in Theorem \ref{thm_canon} the canonical presentation of $v$.
%\end{df}
%
%\begin{comment}
%\begin{df}
%let $\Pi \in A_n$, then the canonical presentation of $\Pi$ is defined to be:
%$$
%\Pi=b_3\cdots b_n, b_i\in B_i
%$$ 
%\end{df}
%\end{comment}

\begin{lem}\label{lem_pres}
%Let $v=(a_1a_2)(a_3a_4)\cdots (a_{k-1}a_k)$, $k \in \bf{N}$ be a permutation in $S_n$ given as a finite product of cycles of length 2 (2-cycles), then:
%\begin{enumerate}
%\item 
%There exists a presentation of $v$ in which $a_1$ appears in the rightmost multuplier (2-cycle) only.
%\item 
%There exists a presentation of $v$ in which $a_2$ appears in the rightmost multuplier (2-cycle) only. 
%\end{enumerate}
Let $k\in \mathbb{N}$, $m_1,\cdots,m_{2k}\in \{1,\cdots,n\}$ be distinct and let $v=(m_1m_2)\cdots(m_{2k-1}m_{2k})\in S_n$, $m_1\neq m_2,\cdots, m_{2k-1}\neq m_{2k}$ be a product of transpositions. Then for every $1\leq i\leq 2k$, there exist a presentation of $v$ as a product of transpositions in which $m_i$ appears in the rightmost factor only. 
\end{lem}

\noindent {\bf Proof of Lemma \ref{lem_pres}}.
We will prove for $i=1$. 
First we will prove for the case $k=2$, i.e. $v$ is a product of two cycles.
If $\{m_1,m_2 \} \cap \{m_3,m_4 \}= \varnothing$ then $v=(m_1m_2)(m_3m_4)=(m_3m_4)(m_1m_2)$ as required. Else, if $m_2=m_3$ then $v=(m_1m_2)(m_2m_4)=(m_1m_2m_4)=(m_2m_4m_1)=(m_2m_4)(m_4m_1)$ as required. The case $m_2=m_4$ is similar. Else, if $m_1=m_3$ then $v=(m_1m_2)(m_1m_4)=(m_2m_1m_4)=(m_4m_2m_1)=(m_4m_2)(m_2m_1)$ as required. The case $m_1=m_4$ is similar. If $m_1=m_3$ and $m_2=m_4$ then $v$ is the identity, thus its cycles are disjoint and commute with each other. All possible cases were checked and thus we finished.   

Now we turn to the general case, where $v$ is a product of $k$ cycles. By induction the Lemma applies also for this case as we can perform the same steps described in the simple case repeatedly until the desired form of $v$ is achieved. \qed

\vskip 0.25 truecm

%\noindent{\bf Proof of Theorem \ref{shuff.1}(1).}
\noindent {\bf Proof of Theorem \ref{thm_canon}}. By induction on $n$. For $n=3$, $A_3=\{(12)(13),(12)(23),e \}=R_3$ and thus the claim holds. Now assume that each $w \in A_{n-1}$, $n\geq 4$, has a unique canonical presentation $w=w_3\cdots w_{n-1}, \quad w_i\in R_i$. We will show that if $v \in A_n$ then $v$ has a unique canonical presentation as well. This actually follows from Lemma \ref{lem_pres}. We will assume that $v \in A_n \setminus A_{n-1}$, otherwise the proof follows immediately from the induction hypothesis. First we apply Lemma \ref{lem_pres} to $v$ to get $n$ in the rightmost factor only. We have $v=g_1\cdots g_k=g_1\cdots g_{k-1}(12)(12)g_k$ with $n$ in $g_k$ only. Now, since $g_1\cdots g_{k-1}(12) \in A_{n-1}$, according to the hypothesis it has a unique canonical presentation, say $w_1\cdots w_{t}$. Thus we have $v=w_1\cdots w_{t}(12)g_k$ and that is unique canonical presentation for $v$, because $(12)g_k$ is unique. \qed

\vskip 0.25 truecm

By Theorem \ref{thm_canon} we conclude
\begin{cor}
The set of A-transpositions, $T(A_n):=\{a_{ij} \mid 1 \le i < j \le n \}$, generates the alternating group on $n$ letters.
\end{cor}

\begin{df}\label{deflenm}
For $v\in A_n$ with the canonical presentation $v=v_3 \cdots v_n$, let 
$$
\hat{\ell}(v)=\#\{i\mid v_i\neq e\}
$$
\end{df}

%\begin{df}
%Denote the generators set for $A_n$ by 
%$$
%G_n=\{a_{i,j}=(12)(ij)\mid 1\le i < j\le n\}
%$$
%\end{df}
%\begin{df}
%$\ell_{T(A_n)}(v)=\min \{k \mid v= a_1 \cdots a_k, \quad a_i\in T(A_n)\}$
%\end{df}

\begin{thm}\label{thm3} For all $v\in A_n$,
$$
\hat{\ell}(v) = \ell_{T(A_n)}(v).
$$
\end{thm} 
In other words the length of the canonical presentation coincides with the natural length with respect to the generating set $T(A_n)$.

\noindent {\bf Proof of Theorem \ref{thm3}}. It suffices to show that if $\hat{\ell}(v)=r$ then $v$ can not be presented as a product of less than $r$ generators. For $n=3$ it was shown that $A_3=R_3$, thus all the elements in $A_3$ are of length $1$, except for the identity $e$ whose length is $0$. For $n > 3$ denote the length of the canonical presentation of $v$ by 
%canonical presentation of $v \in A_n$ by $v_1=a_1\cdots a_m, \quad a_i\in R_i$. Denote also 
$\hat{\ell}(v)=r$. Denote the shortest presentation of $v$  by $v_2=b_1\cdots b_k, \quad b_i\in T(A_n)$. Then $\ell_{T(A_n)}(v)=k$. Now we can apply the corollary of Lemma \ref{lem_pres} described in the proof of Theorem \ref{thm_canon} to turn $v_2$ into a canonical presentation of $v$, say $v_2^{'}$, with $\ell_{T(A_n)}(v_2^{'})=k$. Since $v_1$ and $v_2^{'}$ are two canonical presentations of the same permutation, according to Theorem 4.2 they are actually the same presentation, i.e. $r=k$. \qed  \\ 

In the rest of this paper we will explore the natural length function with respect to $T(A_n)$. For this purpose we will use the equivalence to the length of the canonical expression proved above, as needed. 

%Having proved that the definition of length according to the canonical presentation, $\ell_m(\cdot)$, coincides with the natural definition of length, $\ell_{T(A_n)}(\cdot)$, we will use only the latter one for the rest of this paper.

%\noindent {\bf Proof} of Theorem 4.6. We will show by induction that if $\ell_m(v)=r$ then $v$ can not be presented as a product of less than $r$ generators. For $i=3$ it was shown that $A_3=R_3$ thus all the elements in $A_3$ are of length $1$ except for the identity $e$ which length is $0$. Now assume that the theorem is correct for each $w \in A_{n-1}$. We will show that the theorem is correct for $v\in A_n$. Denote the canonical presentation of $v$ after removing the identity factors by $v_1=v_{i_1}\cdots v_{i_r}$. Then $\ell_m(v)=r$. Denote the shortest presentation of $v$  by $v_2=a_{i_1}\cdots a_{i_k}$. Then $\ell_A(v)=k$. Now we can apply the conclusion from Lemma 4.4 described in the proof of Theorem 4.2 to turn $v_2$ into a canonical presentation of $v$, say $v_2^{'}$. Note that this process does not change the number of cycles in $v_2$, but only the cycles themselves. Thus we have $\ell_A(v_2)=\ell_A(v_2^{'})$.
%Since $v_1$ and $v_2^{'}$ are two canonical presentations of the same permutation, according to Theorem 4.2 they are actually the same presentation, i.e. $r=k$. \qed 

\section{Counting Permutations in $A_n$}
In this section we study the number of permutations in $A_n$ of a given length with respect to $T(A_n)$. A Stirling-type recurrence relation for this statistic is described. 

\begin{df}
Let 
$$
A(n,m) = \{v\in A_n \mid \ell_{T(A_n)}(v)=m \}
$$
and 
$$
a(n,m)=|A(n,m)|
$$
\end{df}

\begin{pro}\label{prop1} For $n\geq 3$,
$$
a(n,1)=a(n-1,1)+n-1
$$
\end{pro}

\noindent{\bf Proof of Proposition \ref{prop1}}. By Definition \ref{def_newgen}, $R_n \setminus \{e\}$ is the subset of generators of $A_n$ that do not belong to $A_{n-1}$; namely $R_n \setminus \{e\}=(T(A_n)\setminus T(A_{n-1}))=\{(12)(nj)\mid 1\leq j <n\}$. These are the generators that involve the new letter $n$. Thus $|R_n|=n$. 
%Since the generators of $A_{n-1}$ are also generators of $A_n$, we can summarize and have $a(n,1)=a(n-1,1)+n-1$. \qed
For every $n$, $A(n,1)=T(A_n)$, and since $T(A_n)=T(A_{n-1})\cup (R_n \setminus \{e\})$, disjoint union, we conclude $a(n,1)=a(n-1,1)+n-1$.  \qed

%By induction on $n$. $A_2$ has no generators and for $A_3$ it was shown in the proof of \ref{thm_canon} that $A_3$ has $2$ generators. And indeed $a(3,1)=a(2,1)+3-1=0+3-1=2$ as required. Assume $a(n-1,1)=a(n-,1)+n-2$. We will show that $a(n,1)=a(n-1,1)+n-1$. 

\begin{thm}\label{thm_groups}
$$
A(n,m) = A(n-1,m-1)\cdot R_n~ \cup ~A(n-1,m),
$$
disjoint union.
\end{thm}

\noindent{\bf Proof of Proposition \ref{thm_groups}}. By two-sided set inclusion. First we will prove that $A(n,m) \supseteq A(n-1,m-1)\cdot R_n \cup A(n-1,m)$. Note that the right hand side of the equation is a disjoint union, according to the $A_n$ canonical presentation properties.

$A(n-1,m) \subseteq A(n,m)$ because a permutation $v$ of length $m$  in $A_{n-1}$ is also of length $m$ in $A_n$. The new generators in $A_n$ can not shorten the length of $v$ because they involve the new letter $n$ which is a fixed point in $v$. 
 
Let $v\in A(n-1,m-1)$, and consider its canonical presentation. Multiply $v$ by $w \in R_n$ from the right side to have the canonical presentation of a permutation $v\cdot w\in A_n$ of length $m-1+1=m$, i.e. $v\cdot w\in A(n,m)$.

We showed that each part of the union on the right hand side of the equation contained in the left hand side, therefore the union itself is also contained in the left side. That proves the first inclusion.
Now we will show that $A(n,m) \subseteq A(n-1,m-1)\cdot R_n \cup A(n-1,m)$. Let $v\in A(n,m)$
\begin{enumerate}
\item
If $n$ is a fixed point in $v$ then $v\in A(n-1,m)$ with the same canonical presentation. 
\item 
Otherwise, $n$ is not a fixed point and therefore the canonical presentation of $v$ is as follows.
$$
v= \underbrace{r_{1}\cdots r_{k-1}}_\text{$\in A(n-1,m-1)$} \cdot \underbrace{r_{n}}_\text{$\in (R_n \setminus \{e\})$}, \quad r_{i}\in R_{i}, \quad 1\leq i \leq n
$$
$r_{n}\in R_n \setminus \{e\}$ because $n$ is not a fixed point and must appear in the presentation. The above canonical presentation of $v$ shows the required inclusion.
\end{enumerate} \qed

From Proposition \ref{prop1} and  Theorem \ref{thm_groups} we can conclude the following relation.
\begin{cor}\label{conc_len} For $1\leq m\leq n-2$,

$$
a(n,m)=a(n-1,m-1)\cdot (n-1) + a(n-1,m).
$$
\end{cor}

\section{Relation between Length and Cycle Number}
In this section we show that the length $\ell_{T(A_n)}(\cdot)$ and the number of cycles $cyc(\cdot)$ are strongly related statistics on $A_n$.

\begin{obs}\label{obs_cyc}
For $n\geq 3$, and $v\in A(n,1)$, $cyc(v)=n-2$
\end{obs}
In other words, the number of cycles in a generator of $A_n$ is $n-2$.

\noindent{\bf Proof of Observation \ref{obs_cyc}}. Every generator $v\in A_n$ is of the form $(12)(in)$, where $1\leq i<n$. If $i=1$, or $i=2$ then $v$ has one cycle of length $3$ and $n-3$ cycles of length $1$ (fixed points). That sums to $n-2$ cycles. Otherwise $i>2$ and then $v$ has two cycles of length $2$ and $n-4$ more cycles of length $1$. That also sums to $n-2$ cycles in $v$. \qed

%\noindent{\bf Proof} of Proposition \ref{prop_cyc} by induction on $n$. For $n=3$ we have $A_3=\{(12)(13),(12)(23),e \}=\{(213)(123),e \}$ and thus the two permutations of length $1$ in $A_3$ have one cycle. And indeed $m(3,1)=3-2=1$ as expected. We will assume that $m(n-1,1)=n-3$ and prove for $n$. The proof is divided into two cases.
%\begin{enumerate}
%\item 
%Let $v\in A(n,1)\cap A(n-1,1)$, namely $v$ is a generator of $A_n$ and of $A_{n-1}$. That makes $n$ a fixed point in $v$ and therefore the cycle $(n)$ is the only addition to the cycle form of $v$ in $A_{n-1}$. Using the induction hypothesis, we get in this case $m(n,1)=m(n-1,1)+1=n-3+1=n-2$ as expected.
%\item 
%Let $v\in R_n$, namely a generator of $A_n$ but not of $A_{n-1}$. Recall the definition $R_n= \{ (12)(nj) \mid 1 \leq j<n \} \cup \{ e \}$. If $j=1$, or $j=2$ then $v$ has one cycle of length $3$ and $n-3$ cycles of length $1$ (fixed points). That sums to $n-2$ cycles. If $j>2$ then $v$ has two cycles of length $2$ and $n-4$ more cycles of length $1$. That also sums to $n-2$ cycles in $v$, as expected.
%\end{enumerate}
%We showed that in both cases, a generator $v\in A_n$ has $n-2$ cycles and therefore proved the claim. \qed

\begin{cor}\label{pro_nice}
For every $n\geq 2$ and $v\in A_n$, the letters $1,2$ are in the same cycle in $v$ if and only if $\ell_{T(A_n)}(v)$ is odd.
\end{cor}

\noindent{\bf Proof of Corollary \ref{pro_nice}}. If $v\in A_n$ is of length one, $1,2$ share the same cycle since the structure of a generator is $(12)(ij)$. In length two, $1,2$ appear in different cycles because of the multiplication process described in the proof of Theorem \ref{thm_groups}. For length three, $1,2$ are in the same cycle according to the same process, and so on and so forth. For odd length, the letters $1,2$ are in the same cycle and for even length they are in a different cycles. That proves both sides of the proposition. \qed

\begin{thm}\label{thm_cyc}
For $n\geq 3$ and $v\in A_n$,
%cyc(v)=n-\ell_{T(A_n)}(v) +\frac{(-1)^{\ell_{T(A_n)}(v)}-1}{2}=\\
$$
cyc(v)=
\begin{cases}
n-\ell_{T(A_n)}(v)& \text{if $\ell_{T(A_n)}(v)$ is even}, \\
n-\ell_{T(A_n)}(v)-1&  \text{if $\ell_{T(A_n)}(v)$ is odd}.
\end{cases} 
$$
\end{thm}

\noindent{\bf Proof of Theorem \ref{thm_cyc}}. By induction on $n$. For $n=3$, $A_3=\{(12)(13)=(213),(12)(23)=(123),e=(1)(2)(3) \}$ and the claim follows. Assume that for each $v\in A_n$
$$
cyc(v)=
\begin{cases}
n-\ell_{T(A_n)}(v)& \text{if $\ell_{T(A_n)}(v)$ is even}, \\
n-\ell_{T(A_n)}(v)-1&  \text{if $\ell_{T(A_n)}(v)$ is odd}.
\end{cases} 
$$
Now, $w\in A_{n+1}$ can be obtained in two ways by Theorem \ref{thm_groups}. First, by multiplying $v\in A_n$ by $r\in R_{n+1}$, and secondly by adding the letter $n+1$ as fixed point to some $v\in A_n$. Both cases will be analyzed.
\begin{enumerate}
\item 
In this case $w=vr$ for $v\in A_n,\quad r\in R_{n+1}$. If $\ell_{T(A_n)}(v)$ is even then, by Corollary \ref{pro_nice}, the letters $1,2$ are in different cycles in $v$ and therefore they will be in the same cycle in $w$, thus $cyc(w)=cyc(v)-1$ (see Theorem \ref{thm_groups} for details). The length of $w$ is $\ell_{T(A_{n+1})}(w)=\ell_{T(A_n)}(v)+1$. By the induction hypothesis,
$$
cyc(w)=cyc(v)-1=n-\ell_{T(A_n)}(v)-1=n-\ell_{T(A_{n+1})}(w)=(n+1)-\ell_{T(A_{n+1})}(w)-1,
$$ 
as required. If $\ell_{T(A_n)}(v)$ is odd then, by Corollary \ref{pro_nice} and Theorem \ref{thm_groups}, $cyc(w)=cyc(v)+1$ and $\ell_{T(A_{n+1})}(w)=\ell_{T(A_n)}(v)+1$.  By the induction hypothesis,
$$
cyc(w)=cyc(v)+1=n-\ell_{T(A_n)}(v)-1+1=n-\ell_{T(A_{n+1})}(w)+1=(n+1)-\ell_{T(A_{n+1})}(w),
$$
as required.
\item 
In this case $w=v$ for some $v\in A_n$, where the letter $n+1$ is a fixed point in $w$. Here, $cyc(w)=cyc(v)+1$ and $\ell_{T(A_{n+1})}(w)=\ell_{T(A_n)}(v)$. If $\ell_{T(A_n)}(v)$ is even,
$$
cyc(w)=cyc(v)+1=n-\ell_{T(A_n)}(v)+1=n-\ell_{T(A_{n+1})}(w)+1=(n+1)-\ell_{T(A_{n+1})}(w),
$$
as required, and if $\ell_{T(A_n)}(v)$ is odd,
$$
cyc(w)=cyc(v)+1=n-\ell_{T(A_n)}(v)-1+1=n-\ell_{T(A_{n+1})}(w)=(n+1)-\ell_{T(A_{n+1})}(w)-1,
$$
as required.
\end{enumerate}
In both cases the relation between cycle number and length holds, therefore the Theorem is proved.\qed
\\

The following relation is By Corollary \ref{pro_nice} and Theorem \ref{thm_cyc}.

\begin{cor}\label{conc_cyc}
Let $v\in A_n$.
\begin{equation}\label{len_cyc}
\ell_{T(A_n)}(v)= \\
\begin{cases}
n-cyc(v)&  \text{if 1,2 are in different cycles of v}, \\
n-cyc(v)-1&  \text{if 1,2 in the same cycle of v}
\end{cases}
\end{equation}
\end{cor}

%\noindent{\bf Proof} of Corollary \ref{conc_cyc}. By Corollary \ref{pro_nice} and Theorem \ref{thm_cyc}.\qed

Equation (\ref{len_cyc}) provides a simple way to find the length of a permutation $v$ given as a product of disjoint cycles. \\

Theorem \ref{thm_cyc} implies that all the permutations of the same length in $A_n$ have the same number of cycles.
\begin{df}
$$
m(n,k) = \emph{number of cycles in a permutation $v\in A_n$ of length $\ell_{T(A_n)}(v)=k$}
$$
\end{df}

\section{Generating Function of Length in $A_n$}\label{secgenfun}
An explicit formula for the generating function of the length in $A_n$, with respect to the generating set $T(A_n)$, is given in this section.

%\begin{obs}\label{obs1}
According to Theorem \ref{thm_cyc}, the number $m(n,k)$, of cycles in a permutation $v\in A_n$ of length $k$, can be calculated by the following formula.
\begin{equation}\label{eq_gen1}
m(n,k)= \\
\begin{cases}
n-k,& \text{if $k$ is even}, \\
n-k-1,&  \text{if $k$ is odd.}
\end{cases} 
\end{equation}
A well known result from $S_n$ is
\begin{equation}\label{eq_gen2}
m(n,k)=n-k
\end{equation}
where the length $k$ is taken with respect to the generating set $T=\{(ij)\mid 1\leq i < j\leq n\}$, namely all the transpositions in $S_n$. Since the number of cycles in a permutation is independent of the generating set, we can conclude from equations (\ref{eq_gen1}) and (\ref{eq_gen2}) that for $v\in A_n$
\begin{equation}\label{eq_gen3}
\ell_{T}(v)=
\begin{cases}
\ell_{T(A_n)}(v),& \text{if $\ell_{T(A_n)}(v)$ is even}, \\ 
\ell_{T(A_n)}(v)+1,& \text{if $\ell_{T(A_n)}(v)$ is odd,}
\end{cases} 
\end{equation}
where $\ell_{T}(v)$ is the length with respect to $T$.
Note that in each of the cases $\ell_{T}(v)$ is even, which complies with the fact that we deal with even permutations in $S_n$.
From equation (\ref{eq_gen3}) we can conclude that the number of permutations  of even length $k$ in $S_n$ equals the sum of the number of permutations of lengths $k$ and $k-1$ in $A_n$. Since the number of permutations of length $k$ in $S_n$ with respect to $T$ is the unsigned Stirling number of the first kind $c(n,n-k)$, we can deduce the following equation for even $k$.
\begin{equation}\label{eq_gen4}
c(n,n-k)=a(n,k)+a(n,k-1).
\end{equation}
%\end{obs}
Furthermore,
\begin{cla}\label{cla1}
Equation (\ref{eq_gen4}) holds also for odd $k\in \mathbb{N}$.
\end{cla}

\noindent{\bf Proof of Claim \ref{cla1}}. Let $k+1$ be even. Using the recursive relation of Stirling numbers 
%(see claim \ref{str1rec}) 
we can develop the left hand of equation \ref{eq_gen4} to have
\begin{align*}
%\lefteqn
&c(n,n-(k+1))=c(n-1,n-(k+1)-1)+c(n-1,n-(k+1))\cdot(n-1)= \\
&c(n-1,n-1-(k+1))+c(n-1,n-1-((k+1)-1))\cdot (n-1)= \\
&c(n-1,n-k-2)+c(n-1,n-1-k)\cdot (n-1)
\end{align*}

Switching sides gives the following result.
$$
(n-1)\cdot c(n-1,n-1-k)=c(n,n-k-1)-c(n-1,n-k-2)
$$
The expressions at right hand side represent even length, so we can use equation \ref{eq_gen4} and conclusion \ref{conc_len} to obtain the desired result.
\begin{align*}
&(n-1)\cdot c(n-1,n-1-(k-1))= a(n,k)+a(n,k-1)-a(n-1,k)-a(n-1,k-1)= \\
&a(n-1,k-1)\cdot (n-1) + a(n-1,k)+a(n-1,k-2)\cdot (n-1) + a(n-1,k-1)-a(n-1,k)-a(n-1,k-1) \\
&=a(n-1,k-1)\cdot (n-1)+a(n-1,k-2)\cdot (n-1)
\end{align*}
Divide both hand sides by $(n-1)$ the following is deduced, for odd $k$.
$$
c(n-1,n-1-k)=a(n-1,k)+ a(n-1,k-1)
$$
\qed

%In the substitution method, a variable is isolated (solved for) and then used to substitute in other equation(s). 
%e.g. 
%2x + y - 6 = 0 
%x + 3y - 13 = 0 
%Solve for y in the first equation. Move the 2x and -6 to the right side, switch signs as we switch sides. 
%y = -2x + 6 
%Now take this definition of y and substitute it in the second equation; any place you see y, make the substitution. 
%x + 3(-2x + 6) - 13 = 0 
%x + -6x + 18 - 13 = 0 
%Collect like terms and simplify. 
%-5x + 5 = 0 
%-5x = -5 
%x = 1 
%Substitute this value of x into the definition for y above. 
%y = -2(1) + 6 
%y = -2 + 6 
%y = 4 

The following generating function for unsigned Stirling numbers of the first kind is well known~\cite [pp. 213]{CO}.
\begin{equation}\label{str_book}
\sum_{k=1}^n c(n,k)\cdot x^{n-k}=(1+x)(1+2x)\cdots (1+(n-1)x)
\end{equation}

By equation (\ref{eq_gen4}),
$$
\sum_{k=0}^n c(n,n-k)\cdot x^{k}=\sum_{k=0}^n a(n,k)\cdot x^{k} + \sum_{k=0}^n a(n,k-1)\cdot x^{k}
$$

Using equation \ref{str_book} for the left hand side, we have the following.
\begin{multline*}
(1+x)(1+2x)\cdots (1+(n-1)x)= \\
 \sum_{k=0}^n a(n,k)\cdot x^{k} + x\cdot \sum_{k=0}^n a(n,k-1)\cdot x^{k-1}= (1+x)\cdot \sum_{k=0}^n a(n,k)\cdot x^{k}
\end{multline*}
Divide both hand sides by $(x+1)$ to get the generating function of length in $A_n$ with respect to the generating set $T(A_n)$.

\begin{thm}\label{thm_genfunc}
\begin{align}\label{eq_genfunc}
\sum_{k=0}^n a(n,k)\cdot x^{k}
&=(1+2x)(1+3x)\cdots (1+(n-1)x) \\
&= \prod_{t=2}^{n-1}(1+tx) \nonumber
\end{align}
\end{thm}

\section{Expectation and Variance}
In this section the expectation and variance of the length function in $A_n$ will be studied.

\begin{df}
Let $A$ be a finite set, and $s:A\rightarrow \mathbb{R}$ a real function. The expectation of $s$ is defined
$$
E[s]:=\frac{1}{|A|} \sum_{a\in A} s(a)
$$
and the variance of $s$ is defined
$$
Var[s]:=E[s^2]-E^2[s]
$$
\end{df}

Given a generating function of $s$, we can use it to calculate these statistics. The following formulas are well-known.

\begin{pro} \label{pro_proba}
Let 
$$
F_s(x) := \sum_{a\in A} x^{s(a)}
$$ 
be the generating function of $s$. Then 
$$
E[s]=\frac{1}{|A|} F^{'}_s(x)\Bigg\vert_{x=1}
$$ 
and
$$
Var[s]=\frac{1}{|A|}\Big[F^{''}_s(x) + F^{'}_s(x) - \frac{1}{|A|}(F^{'}_s(x))^2\Big]\Bigg\vert_{x=1}
$$
\end{pro}

\begin{df}
Recall the definitions of harmonic numbers (see Definitions \ref{def_harmc} and \ref{def_genharmc}).
$$
H_n= 1+\frac{1}{2}+\frac{1}{3}+\cdots + \frac{1}{n}
$$
$$
H_{n,m}= 1+\frac{1}{2^m}+\frac{1}{3^m}+\cdots + \frac{1}{n^m}
$$
\end{df}

\begin{thm}\label{thm_exp_var}
%Let
%\begin{align}\label{eq_genfunc}
%g_n(x)=\sum_{v\in A_n} x^{\ell_{T(A_n)}(v)}
%= \prod_{t=2}^{n-1}(1+tx)  
%=(1+2x)(1+3x)\cdots (1+(n-1)x) \nonumber
%\end{align}
%be the generating function of length in $A_n$ with respect to the generating set $T(A_n)$ (see Theorem \ref{thm_genfunc}), then
The expected value of $\ell_{T(A_n)}$ is
$$
E[\ell_{T(A_n)}]=n-H_n-\frac{1}{2}
$$
and its variance is
$$
Var[\ell_{T(A_n)}]=H_n-H_{n,2}-\frac{1}{4}
$$
\end{thm}
\noindent{\bf Proof} of Theorem \ref{thm_exp_var}. Compute the derivative of the generating function of length (see Theorem \ref{thm_genfunc}) as a product of functions.
$$
 \bigg(\prod_{t=2}^{n-1}(1+tx)\bigg)^{'}=  \bigg(\prod_{t=2}^{n-1}(1+tx)\bigg)\sum_{t=2}^{n-1}{\frac{t}{1+tx}}.
$$
%Now $E[\ell_{T(A_n)}]$ can be calculated.
Thus, by Proposition \ref{pro_proba},
$$
E[\ell_{T(A_n)}]=\frac{1}{|A_n|}\bigg(\prod_{t=2}^{n-1}(1+tx)\bigg)\sum_{t=2}^{n-1}{\frac{t}{1+tx}}\Bigg\vert_{x=1}=\frac{2}{n!}\cdot \frac{n!}{2}\bigg(n-2-\sum_{t=2}^{n-1}\frac{1}{1+t}\bigg)=n-H_n-\frac{1}{2}
$$
The variance calculation follows.
\begin{align*}
\bigg(\prod_{t=2}^{n-1}(1+tx)\bigg)^{''}
&=\bigg[\bigg(\prod_{t=2}^{n-1}(1+tx)\bigg)\sum_{t=2}^{n-1}{\frac{t}{1+tx}}\bigg]^{'} \\
&=\bigg(\prod_{t=2}^{n-1}(1+tx)\bigg)\sum_{t=2}^{n-1}{\frac{t}{1+tx}}\sum_{t=2}^{n-1}\frac{t}{1+tx} +\bigg(\prod_{t=2}^{n-1}(1+tx)\bigg)\sum_{t=2}^{n-1}\frac{-t^2}{(1+tx)^2} 
\end{align*}
\begin{align*}
Var[\ell_{T(A_n)}]&=\frac{1}{|A_n|}\bigg[\bigg(\prod_{t=2}^{n-1}(1+tx)\bigg)^{''}+\bigg(\prod_{t=2}^{n-1}(1+tx)\bigg)^{'}-\frac{1}{|A_n|}\bigg(\bigg(\prod_{t=2}^{n-1}(1+tx)\bigg)^{'}\bigg)^2\bigg]_{x=1} \\
&=\frac{2}{n!}\bigg[\bigg(\prod_{t=2}^{n-1}(1+tx)\bigg)\sum_{t=2}^{n-1}{\frac{t}{1+tx}}\sum_{t=2}^{n-1}\frac{t}{1+tx} 
+\bigg(\prod_{t=2}^{n-1}(1+tx)\bigg)\sum_{t=2}^{n-1}\frac{-t^2}{(1+tx)^2} \\
&+\bigg(\prod_{t=2}^{n-1}(1+tx)\bigg)\sum_{t=2}^{n-1}\frac{t}{1+tx} -\frac{2}{n!}\bigg(\bigg(\prod_{t=2}^{n-1}(1+tx)\bigg)\sum_{t=2}^{n-1}\frac{t}{1+tx}\bigg)^2\bigg]_{x=1} \\
&=\frac{2}{n!}\bigg[\frac{n!}{2}\bigg(n-H_n-\frac{1}{2}\bigg)^2 + \frac{n!}{2}\bigg(2H_n+\frac{1}{4}-n-H_{n,2}\bigg)+\frac{n!}{2}(n-H_n-\frac{1}{2}) \\
&-\frac{2}{n!}\bigg(\frac{n!}{2}(n-H_n-\frac{1}{2}\bigg)^2\bigg] \\
&= H_n-H_{n,2}-\frac{1}{4}
\end{align*}
\qed

\section{Connection with Restricted Stirling Numbers}
This section discusses the relation between our statistic $a(n,m)$ and $2$-restricted Stirling numbers of the first kind (see Broder~\cite[\S 1]{BRO}). This relation was initially established using the On-Line Encyclopedia of Integer Sequences~\cite{EIS}.

%\begin{df}
%The unsigned 2-restricted Stirling number of the first kind, denoted ${\ribua{n}{k}}_2$, count the number of permutations on the letters $\{1,2,...,n\}$ of $k$ disjoint cycles, such that letters $1$ and $2$ belong to distinct cycles. This is the particular case $r=2$ of the unsigned r-restricted Stirling numbers of the first kind, defined above in \ref{dfrstr1}.
%\end{df}

\medskip
Recall Definition \ref{dfrstr1} of the $r$-restricted Stirling numbers of the first kind, $\ribua{n}{k}_r$. We shall use it with $r=2$.  $\ribua{n}{k}_1=c(n,k)$ is the usual(unrestricted) Stirling numbers of the first kind.

\begin{cla}\label{cla_str} (See Broder~\cite[\S 3, Thm. 3]{BRO} for a generalized version)
$$
c(n,k)={\ribua{n}{k}}_2+{\ribua{n}{k+1}}_2
$$
\end{cla}

%\noindent{\bf Proof of Claim \ref{cla_str}}. Consider the number $c(n,k)$ which denotes the number of permutations is $S_n$ with exactly $k$ cycles. It can be naturally divided into two numbers. The first counts the permutations where $1,2$ appear in distinct cycles, and the second counts the permutations where $1,2$ appear in the same cycle. The first number is exactly ${\ribua{n}{k}}_2$. We will prove that the second number is ${\ribua{n}{k+1}}_2$ through bijection. For the first direction, Let $v\in A_n$ with $1,2$ in the same cycle $c$. Arrange this cycle so that $1$ will be the leftmost element. The result would have the following pattern.  
%$$
%c=(1,t_1,t_2,\cdots ,t_l,2,k_1,k_2,\cdots ,k_r) \qquad 0\leq r,l\leq n-2
%$$
%Now break this cycle into two cycles. The first is $c_1=(1,t_1,t_2\cdots t_l)$ and the second is $c_2=(2,k_1,k_2\cdots k_r)$. Replace the cycle $c$ in $v$ with $c_1$,$c_2$ to get a permutation $v^{'}\in {\ribua{n}{k+1}}_2$. The described process is deterministic so we have injection to ${\ribua{n}{k+1}}_2$. For the second direction of the bijection do the reverse process. Take the two distinct cycles where $1,2$ appear, and stick them into one cycle to get a permutation $v\in c(n,k)$ where $1,2$ appear in the same cycle. That is also an injection and therefore we proved the bijection. This completes the proof of the claim. \qed

\begin{thm}\label{thm_str}
The number of permutations in $A_n$ of length $\ell_{T(A_n)}(\cdot)=k$ is equal to a corresponding 2-restricted stirling number. Namely,
$$
a(n,k)={\ribua{n}{n-k}}_2, \quad (0\leq k\leq n-2)
$$
\end{thm}

We will give two proofs to Theorem \ref{thm_str}. The first proof is algebraic and the second is a direct bijection between two sets.

\noindent{\bf Proof of Theorem \ref{thm_str}}. From claims \ref{cla_str} and \ref{cla1} we can deduce the following equation.
%$$
%a(n,n-k)+a(n,n-k-1)={\ribua{n}{k}}_2+{\ribua{n}{k+1}}_2
%$$
%For convinience we will replace $k$ by $n-k$.
%Replace $k$ by $n-k$ for convinience:
$$
a(n,k)+a(n,k-1)={\ribua{n}{n-k}}_2+{\ribua{n}{n-k+1}}_2, \quad (0\leq k\leq n-1)
$$ 
Now the Theorem can be proved by induction on $k$. By assumption, $n\geq 2$. For $k=0$ we have $a(n,0)=1$ and ${\ribua{n}{n}}_2=1$. The claim $a(n,k)={\ribua{n}{n-k}}_2$ now follows by induction on $k$.
%$$
%
%$$ 
%For $k$ we have 
%$$
%a(n,k)+a(n,k-1)={\ribua{n}{n-k}}_2+{\ribua{n}{n-k+1}}_2,
%$$ 
%Hence, by the induction hypothesis,
%$$
%a(n,k)={\ribua{n}{n-k}}_2, \quad n\geq 2,
%$$
%Completing the proof. \qed

\begin{df}
Let
$$
P(n,k)=\{v\in S_n \mid cyc(v)=k \text{ and } 1,2 \text{ are in different cycles in } v\}
$$
\end{df}

\medskip
An explicit bijection between the sets $A(n,k)$ and $P(n,n-k)$ will be presented.
\\ \\
\noindent{\bf A Bijective Proof of Theorem \ref{thm_str}}.
Define a map $f:A(n,k)\rightarrow P(n,n-k)$
$$
f(v)=
\begin{cases}
v& \text{if $\ell_{T(A_n)}(v)$ is even} \\
(1,2)v&  \text{if $\ell_{T(A_n)}(v)$ is odd}
\end{cases} 
$$
We will show that $f$ is one-to-one and onto $P(n,n-k)$.
\begin{enumerate}
\item
Consider $v_1,v_2\in A(n,k)$ with $f(v_1)=f(v_2)$. If $v1,v2$ are both of even length, or both of odd length then, by the definition of $f$, $v_1=v_2$.
If $v_1$ is of even length and $v_2$ is of odd length or vice versa then, by the definition of $f$ and the fact that $f(v_1)=f(v_2)$, $v_1=(1,2)v_2$. This contradicts the assumption that $v_1,v_2\in A(n,k)$, and therefore impossible. Since only the first case is feasible, $v_1=v_2$ and $f$ is one-to-one.
\item 
Consider $w\in P(n,n-k)$. The length of $w$ in $S_n$, $\ell_{T}(w)$, is $k$. If $k$ is even then $w\in A_n$. By Corollaries \ref{pro_nice} and \ref{conc_cyc} $\ell_{T(A_n)}(w)=k$, therefore $w\in A(n,k)$ and $f(w)=w$. If $k$ is odd then $(1,2)w\in A_n$. By Corollaries \ref{pro_nice} and \ref{conc_cyc} $\ell_{T(A_n)}((1,2)w)=k$, therefore $(1,2)w\in A(n,k)$ and $f((1,2)w)=w$. This proves that $f$ is onto $P(n,n-k)$. \qed

\end{enumerate}

\section*{Acknowledgements}
\noindent I wish to thank my supervisors Prof.\ Ron M.\ Adin and Prof.\ Yuval Roichman for their professional guidance, patience and for introducing me to the fascinating world of combinatorial research. 
\\ \\
I want to thank my parents for their love and support.
\\ \\
I would like to thank all those who helped me and supported me during this work.

\end{document}